\documentclass[10pt,a4paper]{article}
\usepackage[utf8]{inputenc}
\usepackage[T1]{fontenc}
\usepackage{hyperref}
\usepackage[nottoc]{tocbibind} 
\usepackage{fancyhdr}
\usepackage{graphicx}
\usepackage{color}
\usepackage{bm} 
\usepackage{amsmath}
\usepackage{amssymb} 
\usepackage{amsthm}  
\usepackage{mathrsfs} 
\usepackage{ifthen}
\usepackage{tikz}
\usepackage[framemethod=TikZ]{mdframed} 

\newtheorem{definition}{Definition}

\newtheorem{theorem}{Theorem}
\newtheorem{proposition}{Proposition}

\newcommand{\w}[1]{\pmb{#1}}  

\newcommand{\B}[1]{\overline{#1}}
\newcommand{\ts}{\otimes}

\newcommand{\Div}{\w \nabla \cdot}

\newcommand{\la}{\langle}
\newcommand{\ra}{\rangle}

\newcommand{\wnab}{\w{\nabla}}

\newcommand{\dd}{\mathbf{d}}

\newcommand{\defin}[1]{\textbf{\itshape #1}}
\newcommand{\be}{\begin{equation}}
\newcommand{\ee}{\end{equation}}
\newcommand{\bea}{\begin{eqnarray}}
\newcommand{\eea}{\end{eqnarray}}
\newcommand{\bsea}{\begin{subeqnarray}}
\newcommand{\esea}{\end{subeqnarray}}

\newcommand{\RN}[1]{\textup{\uppercase\expandafter{\romannumeral#1}}}

\newcommand{\R}{\mathbb{R}}
\newcommand{\Z}{\mathbb{Z}}

\newcommand{\M}{\mathscr{M}}

\newcommand{\g}{\w{g}}

\newcommand{\T}{\mathrm{T}}

\newcommand{\spose}[1]{\hbox to 0pt{#1\hss}}
\newcommand{\lta}{\mathrel{\spose{\lower 3pt\hbox{$\mathchar"218$}}
     \raise 2.0pt\hbox{$\mathchar"13C$}}}
\newcommand{\gta}{\mathrel{\spose{\lower 3pt\hbox{$\mathchar"218$}}
     \raise 2.0pt\hbox{$\mathchar"13E$}}}

\newcommand{\N}{\mathscr{N}}
\newcommand{\C}{\mathscr{C}}
\newcommand{\U}{\mathscr{U}}
\newcommand{\F}{\mathscr{F}}
\newcommand{\V}{\mathscr{V}}
\newcommand{\W}{\mathscr{W}}

\newcommand{\Hor}{\mathscr{H}}

{\begin{svgraybox}\begin{example}}%
{\end{example}\end{svgraybox}}

 {\small\begin{remark}}%
 {\end{remark}\normalsize}

\newmdenv[backgroundcolor=gray!20!white,roundcorner=5pt,hidealllines=true]{greybox}

\begin{document}

\title{Multi-Fibers Bundles as a new model for high-dimensional Spacetimes.}
\author{Stephane Collion and Michel Vaugon.\\ \sl \small Institut de Math{\'e}matiques, Universit{\'e} Paris VI, Equipe
  G{\'e}om{\'e}trie et Dynamique, 
    \\ \sl \small email: vaugon.michel@wanadoo.fr, stephane.collion@wanadoo.fr}
\date{August 2017, August 2022}
\maketitle

\begin{abstract}
Around 1920, Kaluza and Klein had the idea to add a fifth dimension to the classical 4-dimensional spacetime of general relativity to create a geometric theory of gravitation and electromagnetism. 
Today, theoretical evidences, like string theory, suggest the need for a spacetime with more than five dimensions. 
The mathematical translation of the heuristic idea of a 4-dimensional classical spacetime equipped with extra "small" dimensions, is a fiber bundle structure $\pi : \M\rightarrow \B\M$ on a $(4+k)$-dimensional manifold $\M$, with fiber a compact manifold $\F$ of dimension $k$, more shortly a $\F$-fibration.

Kaluza and Klein used a fibration with fiber the standard circle $S^1$, this fiber carrying the electromagnetic potential. Inclusion of other physical interaction would therefore require a $\F$-fibration with $\F$ of the form $\F=S^1\times W$, $W$ being a compact manifold.

However, with such a $\F$-fibration, they is no naturally defined fiber diffeomorphic to $S^1$ at each point of the manifold, and therefore one looses the possibility to define simply the electromagnetic potential. 

We want to present in this paper a mathematical structure generalizing the fiber bundle, that enable the possible 
definition of multiple naturally defined fibers at each point of the manifold, on which therefore one can define objects that depend only on one of the components, $S^1$ or $W$, of the global $(S^1\times W)$ fiber.

Although we do not pretend here to model precisely other known physical interactions, we present 
this geometric structure as a possible way to model or encode deviations from standard 4-dimensional General Relativity, or "dark" effects such as dark matter or energy ;  (we refer to the authors' article [3] from which this paper is extracted, but whose purpose is different). Also this geometry was a starting point for the second author's  new approach to a geometric unification of General Relativity and Quantum Physics ( see [19]). 
\footnote{Version 02 08 2022. AMS subject classification: 83C22, 83E05, 83E15}\\

\end{abstract}

\footnotesize{ \emph{In all the paper,  $(\M,\g)$ is  a semi-Riemannian manifold with metric $\g$, $\wnab$ is its Levi-Civita connection, $\w R$ and $\w{Ric}$ are the Riemann and Ricci curvatures. $\T_p\M$ is the tangent space of $\M$ at $p$ ; $T_pf$ is the tangent map at $p$ of a function $f$. We sometimes note the scalar product $\g(X,Y):=\la X,Y \ra$. We note $\Div \w T$ the divergence of a tensor $\w T$. In all the paper, using musical isomorphisms, we identify without any comments a $(0,2)$-tensor  $\w T$ with the $(2,0)$-tensor $\w T^{\sharp \sharp}$. They represent the same physical object. In particular, the Einstein curvature $\w G=\w{Ric}-1/2\w S.\g$ will often be consider as a $(2,0)$-tensor, $\w G=G^{ij}$. We also consider, for a 2-tensor $\w T$, the divergence $\Div \w T$ as a vector, that is, we identify the 1-form $\Div \w T$ and $(\Div \w T)^\sharp$. At last,  we shall note $^e \w T$ the endomorphism field $\g$-associated to $\w T$ : 
 $\forall (\w{u},\w{v}) \in \T_p(\M)\times\T_p(\M)\,$ :  $\, \g ( \w{u} , ^e\w{T}(\w{v}))  =\w{T}(\w{u},\w{v})$.}}
 \normalsize

\section{Beyond five dimensions.}
Theoretical evidences, like string theory, suggest the need for a spacetime with more than five dimensions. We want to present in this section a possible extension of our model, that preserve the results obtained so far for the inclusion of electromagnetism, but that enable the possible inclusion of such other dimensions that might model geometrically other physical effects.
Although we do not pretend here to model precisely other known physical interactions, we present a geometric structure giving a possible way, for instance, to model or encode deviations from standard 4-dimensional General Relativity, or "dark" effects such as dark matter or energy.

\subsection{Multi-fiber bundle.}
\subsubsection{Multiple fibers at each point of spacetime.}
The mathematical translation of the heuristic idea of a 4-dimensional classical spacetime equipped with extra "small" dimensions is a fiber bundle structure $\pi : \M\rightarrow \B\M$ on a $(4+k)$-dimensional manifold $\M$, with fiber a compact manifold $\F$ of dimension $k$, more shortly a $\F$-fibration.

Now if we want to keep the result obtained for electromagnetism while including other possible interactions, the fiber $\F$ should be of the form $\F=S^1\times W$ where  $W$ is a compact manifold of dimension $m$ and  $S^1$ is the classical circle.  However, if we want to keep results obtained for electromagnetism in 5 dimensions, through objects naturally given by the action of $S^1$, we face a important issue : at each point $x\in \M$, such a fiber bundle gives a natural fiber $\F=S^1\times W$ through $x$, but it does not give a natural fiber through $x$ isomorphic to $S^1$ only ; there is no natural splitting of the fiber $S^1\times W$  at $x$. Therefore, such a fiber bundle alone will not furnish an electromagnetic potential $\w Y$. 

A very elegant extension of the structure of fiber bundle, giving a way to define any number of natural fibers at each point $x$ of a manifold $\M$, was originally proposed by Michel Vaugon in [ref]. We give here a new approach, based on the more classical notions of  fibrations and submersions.

\subsubsection{Splitting of a product manifold.}
A natural way to split a fiber $\F$ of the form $S\times W$ is based on the following nice construction. Let $\F$, $S$ and $W$ be three compact manifolds and let 
$$\Phi=(h,f) : \F\rightarrow S\times W$$
 be a diffeomorphism, where $h$ and $f$ are the components of $\Phi$. Then $h :\F\rightarrow S$ and $f : \F \rightarrow W$ are submersions. 
We then define unambiguously, for any $x\in\F$, two fibers at $x$ by :
$$ S_x:= f^{-1}(f(x))$$
$$ W_x:=h^{-1}(h(x))$$
Because $S$ and $W$ are compact, a theorem of Ehresmann states that these submersions are in fact fibrations.
Using the diffeomorphism $\Phi$, we can see that $h$ and $f$ are more precisely fibrations with fibers $W$ and $S$ respectively. Indeed, the restrictions 
$h|_{S_x} : S_x\rightarrow S$ and $f|_{W_x} : W_x\rightarrow W$ are diffeomorphisms whose inverse maps are given respectively by :
$$(h|_{S_x})^{-1}(u)=\Phi^{-1}(u,f(x)). \quad \quad  (f|_{W_x})^{-1}(v)=\Phi^{-1}(h(x),v).$$
Thus :
\begin{itemize}
\item $h :\F \rightarrow S$ is a $W$-fibration.
\item $f :\F \rightarrow W$ is a $S$-fibration.
\end{itemize}
Furthermore, we have a natural splitting of the manifold $\F$ as a product of two fibers at a given point : for a given point $p \in\F$, we have a natural diffeomorphism:
\begin{align}
\psi_p : \quad & \F \longrightarrow \quad   S_p \times W_p \nonumber \\
 & y \longmapsto \quad  (\, \Phi^{-1}(h(y),f(p))\, , \, \Phi^{-1}(h(p),f(y))\, ) \nonumber
\end{align}
Indeed, the inverse map is given by :
$$\psi_p^{-1}: S_p \times W_p \longrightarrow  \F$$
$$\quad \quad \quad \quad \quad \quad (a,b) \longmapsto \Phi^{-1}(h(a)\, , \, f(b))$$
To prove this, note first that, for $(a,b)\in S_p \times W_p$, by definition of the fibers, $f(a)=f(p)$ and $h(b)=h(p)$. Setting $y=\Phi^{-1}(h(a) , f(b))$, by definition of $\Phi$ and its components $h$ and $f$, $h(y)=h(a)$ and $f(y)=f(b)$. Therefore, $\Phi^{-1}(h(y), f(p))= \Phi^{-1}(h(a),  f(a))=a$ and $\Phi^{-1}(h(p), f(y))= \Phi^{-1}(h(b),  f(b))=b$.

\subsubsection{Multi-fiber bundle.}
We now apply this construction to define a \defin{multi-fiber} structure on a manifold $\M$ equipped with a classical fiber-bundle structure $\pi : \M\rightarrow \B\M$  whose fiber $\F=S\times W$ is a product of two compact manifolds. By definition, for any point $p\in \M$, there exist a \emph{bundle chart} $(\B \U, \phi)$ with $\pi(p)\in \B \U \subset \B\M$ :
\begin{align}
& \pi^{-1}(\B \U)   \stackrel{\phi}{\longrightarrow}\B \U\times \F =\B\U\times S \times W  \nonumber \\
\pi & \quad \downarrow \quad\quad  \swarrow p_1  \nonumber \\
& \quad \B \U \nonumber
\end{align}

Because of the commutativity of the diagram, $\phi$ can be written $\phi=(\pi,\Phi)$, with $\Phi: \pi^{-1}(\B \U) \rightarrow \F$ and $\Phi=(h,f)$, where $h : \pi^{-1}(\B \U) \rightarrow S$ and $f : \pi^{-1}(\B \U) \rightarrow W$ are the components of $\Phi$ as in the previous section. (To be rigorous, $\phi=(\pi |_{\pi^{-1}(\B\U )},\Phi)$. We also have here the following commutative diagram : 
\begin{align}
& \quad\pi^{-1}(\B \U)    \nonumber \\
(\pi,h)\swarrow \quad &  \quad \downarrow (\pi,\Phi) \quad\quad  \searrow (\pi,f)  \nonumber \\
\B\U \times S \quad\quad \longleftarrow \,\,\, &\B \U\times S\times W  \longrightarrow\quad\quad\B\U \times W \nonumber
\end{align}
where the horizontal arrows on the last line are the obvious projections.

We could therefore, using  $\Phi|_{\F_{ p}} : \F_{ p} \stackrel{\simeq}{ \rightarrow} S\times W$, where $\F_{ p}:=\pi^{-1}(\pi(p))$,  define fibers $S_p$ and $W_p$  :
 $ S_p := (f|_{\F_{p}})^{-1}(f(p))$ and $ W_p := (h|_{\F_{p}})^{-1}(h(p))$.

However, another chart $(\B \U', \phi')$ around $\B p:=\pi(p)$, with $\phi'=(\pi,\Phi')$ and $\Phi'=(h',f')$,  could give rise to different fibers $S_p$ and $W_p$ if for example $(f|_{\F_{p}})^{-1}(f(p))\neq(f'|_{\F_{p}})^{-1}(f'(p))$.

To get well-defined fibers $S$ and $W$ through $p$, we therefore need to impose a compatibility condition between the charts. This will lead to our definition of \defin{multi-fiber bundle}.

To get there, remember the following facts : as a fiber bundle, $\pi : \M\rightarrow \B\M$ can be considered as being equipped with a \emph{bundle atlas}, that is, a family $\{(\B\U_\alpha, \phi_\alpha)\}_{\alpha \in A}$ of bundle charts such that $\{\B\U_\alpha\}_{\alpha \in A}$ is a cover of $\B \M$. Then, if $\B\U_\alpha \cap \B\U_\beta$ is not empty, we have an \emph{overlap map} 
$$\phi_\alpha \circ\phi_\beta^{-1} : (\B\U_\alpha \cap \B\U_\beta)\times\F \rightarrow (\B\U_\alpha \cap \B\U_\beta)\times\F.$$
Also, still writing $\phi_\alpha=(\pi, \Phi_\alpha)$, $\Phi_\alpha |_{\F_p}$  is a diffeomorphism for each $p$ with $\pi(p)\in\B\U_\alpha$. Therefore 
$\Phi_\alpha |_{\F_p}\circ \Phi_\beta |_{\F_p}^{-1} :\F\rightarrow \F$ is a diffeomorphism for all $p$ such that $\pi (p) \in \B\U_\alpha \cap \B\U_\beta$.

It is on these overlap maps that we shall impose a compatibility condition : 
\begin{greybox}
Let $S$ and $W$ be two compact manifolds, and $\F:=S\times W$.

Let $\pi : \M\rightarrow \B\M$ be a fiber-bundle with fiber $\F$, and let $\{(\B\U_\alpha, \phi_\alpha)\}_{\alpha \in A'}$ be a complete bundle atlas for $\pi$. We have at each $p\in\M$ a global $\pi$-fiber : $\F_{ p}:=\pi^{-1}(\pi(p))$

Let's write $\phi_\gamma=(\pi,\Phi_\gamma)$ and $\Phi_\gamma=(h_\gamma, f_\gamma)$ for any bundle chart. We say that $\M$ is a \defin{multi-fiber bundle} with fibers $(S,W)$, or a $(S,W)$-fibration, if there exists a \emph{sub-atlas} $\{(\B\U_\alpha, \phi_\alpha)\}_{\alpha \in A}$ such that we have, for any $\alpha, \beta \in A$ with $\B\U_\alpha \cap \B\U_\beta\neq \emptyset$ and any $p\in \pi^{-1}(\B\U_\alpha \cap \B\U_\beta)$ :
\begin{itemize}
	\item $\phi_\alpha^{-1}(\{\pi(p)\}\times S\times \{f_\alpha(p)\})=\phi_\beta^{-1}(\{\pi(p)\}\times S\times \{f_\beta(p)\})$
	\item $\phi_\alpha^{-1}(\{\pi(p)\}\times\{h_\alpha(p)\}\times W)=\phi_\beta^{-1}(\{\pi(p)\}\times\{h_\beta(p)\}\times W)$
\end{itemize}
It is fair to call this sub-atlas $\{(\B\U_\alpha, \phi_\alpha)\}_{\alpha \in A}$ a multi-fiber atlas, and its charts, muti-fiber charts. With it, we can define unambiguously two new fibers at each point $p\in \M$ : using any chart $(\B\U_\alpha, \phi_\alpha)$ of this multi-fiber atlas with $\pi(p)\in\B\U_\alpha$ : 
\begin{itemize}
\item The $S$-fiber $ S_p := \phi_\alpha^{-1}(\{\pi(p)\}\times S\times \{f_\alpha(p)\})$, a submanifold of  $\F_{p}$
\item The $W$-fiber $\phi_\alpha^{-1}(\{\pi(p)\}\times\{h_\alpha(p)\}\times W)$, a submanifold of  $\F_{p}$
\end{itemize}
\end{greybox}

\emph{Notation : We shall  note $\B x=\pi(x)$ for $x$ in $\M$, and 
$\F_{\B x}:=(S \times W)_{\B x}:=\pi^{-1}(\pi(x)):= \pi^{-1} (\B x)$
the $\pi$-fiber at a point $x$ of $\M$. $(S \times W)_{\B x}$ is intuitively $\{\B x\}\times \F$. We also note $\B\U:=\pi(\U)$ for a set $\U$. For a map $f:\M\rightarrow \N$, where $\N$ is a manifold, we sometimes note $f_{\B x}$ the restriction of $f$ to the $\pi$-fiber $(S \times W)_{\B x}$ : $f_{\B x}:=f|_{\pi^{-1}(\B x)} :=f|_{\pi^{-1}(\pi(x))} :=f|_{\F_{\B p}}$.}

Note that the compatibility condition given in the definition of a multi-fiber bundle is equivalent to the following on the $\Phi_\alpha=(h_\alpha, f_\alpha)$'s :
\begin{greybox}
For any $p$ with $\pi (p) \in \B\U_\alpha \cap \B\U_\beta$, we have :
\begin{itemize}
	\item $\Phi_\alpha |_{\F_p}^{-1} (S\times \{f_\alpha(p)\})=\Phi_\beta |_{\F_p}^{-1} (S\times \{f_\beta(p)\})$
	\item $\Phi_\alpha |_{\F_p}^{-1} (\{h_\alpha(p)\}\times W)=\Phi_\beta |_{\F_p}^{-1} (\{h_\beta(p)\}\times W)$
\end{itemize}	
\end{greybox}

Using this and the splitting of each fiber $\F_{p}$ with the diffeomorphism $\Phi|_{\F_{ p}} : \F_{ p} \stackrel{\simeq}{ \rightarrow} S\times W$ of any multi-fiber chart, as seen in the previous section, we can give equivalent characterization of the fibers $S_p$ and $W_p$ for $p\in\M$ :

\begin{greybox}
\begin{itemize}
\item $ S_p := \Phi_\alpha |_{\F_p}^{-1} (S\times \{f_\alpha(p)\})=(f_\alpha|_{\F_{ p}})^{-1}(f_\alpha(p))$
\item $ W_p := \Phi_\alpha |_{\F_p}^{-1} (\{h_\alpha(p)\}\times W) = (h_\alpha|_{\F_{p}})^{-1}(h_\alpha(p))$
\end{itemize}
\end{greybox}

\begin{greybox}
The fundamental idea of multi-fiber structure is that, for a $(S\times W)$-fibration structure on a manifold $\M$, one can furthermore  define unambiguously, at each point, objects that depend only on one of the components, $S$ or $W$, of the global $(S\times W)$ fiber. See below.
\end{greybox}

The multi-fiber structure satisfies the following natural and important properties :
\begin{greybox}
\begin{itemize}
\item \emph{Fibers are well defined : } $p'\in S_p \Rightarrow S_{p'}=S_p$ and $p'\in W_p \Rightarrow W_{p'}=W_p$.
\item \emph{Splitting of the fibers : }Thanks to the splitting $\psi_p$ defined in the previous section,  we have at each point $p\in \M$ and for each multi-fiber chart $(\B\U_\alpha, (\pi,\Phi_\alpha))$ around $\B p$, a canonical isomorphism  :
$\psi_{p,\alpha} : \F_{p} \stackrel{\simeq}{ \rightarrow} S_p \times W_p$. The choice of another chart $(\B\U_\beta, (\pi,\Phi_\beta))$ around $\B p$ will lead to the same splitting but, via the overlap map, through diffeomorphisms of $S_p$ and $W_p$.
\item \emph{Adapted charts :} Let us fix a point $p\in\M$. A multi-fiber chart around $\B p$ gives a diffeomorphism 
$$\phi_{\U} :\U\rightarrow \B\U \times \F_{p}$$
 from a neighborhood $\U=\U_p$ of $p$ in $\M$. Composing with $Id\times \psi_p$, where $\psi_p$ is the splitting of the fiber  $ \F_{p}$ associated to $\Phi_\U$, we get an adapted diffeomorphism 
$$\varphi_{\U} =(Id\times \psi_p)\circ  \phi_{\U} : \U \rightarrow \B\U \times  S_p \times W_p$$
singularizing the fibers at $p$.  Then, taking coordinates charts on open subsets of  $ \B\U$, $S_p$ and $W_p$ respectively, we obtain very useful adapted charts on $\U_p$. See below.
\item \emph{Orientation of the fibers :} Let  $S$ be given an orientation ; an orientation on $W$ would be treated the same way. For any point $p\in\M$, the diffeomorphisms $h_{\alpha}|_{S_p} : S_p\rightarrow S$ can be used to pull-back the orientation of $S$ on $S_p$. We say that the multi-fiber structure  is compatible with the orientation of $S$ if, in some neighborhood $\U_p$ of any point $p$, there is a frame field for $\T_y S_y$, $y\in\U_p$, compatible with the orientation pulled-back by any $h_{\alpha}|_{S_y}$ such that $y\in\U_{\alpha}\cap\U_p$. 
\end{itemize}
\end{greybox}

We now consider $\M$ to be equipped with a metric $\g$. We then define the \defin{horizontal space} $H_p$ at a point $p\in\M$ as the $\g$-orthogonal space to $\T_p \F_{p}$ in $\T_p\M$ :
\begin{greybox}
$$H_p := (\T_p \F_{p})^\bot $$
\end{greybox}
We can then consider a compatibility condition between $\g$ and the multi-fiber structure:
\begin{greybox}
\begin{itemize}
\item \emph{Signature of the fibers :} We say that the metric  $\g$ is compatible with the multi-fiber bundle structure if the signature of the restriction of $\g$ to any fiber as defined above is constant. That is, for any $p \in M$, the signature of $\g$ restricted to $S_p$ and $W_p$ is independent of $p$. In this case, for given signatures $\sigma_a$ and $\sigma_b$ of adequate length, we will say that $\g$ is of signature $\sigma_a$ on $S$ and $\sigma_b$ on $W$. The signature of $\g_p$ on the horizontal space $H_p$ is then also independent of $p$.
\end{itemize}
\end{greybox}

Note that all the above construction can easily be generalized to define more than 2 fibers at each point of a manifold $\M$. If we are given a global $\pi$-fiber of the form $\F=W_1\times ... \times W_k$, we essentially replace the compatibility conditions on the charts by something like :  
$$\phi_\alpha^{-1}(\{\pi(p)\}\times\{f^1_\alpha(p)\}\times...\times W_j\times...\times \{f^k_\alpha(p)\})=\phi_\beta^{-1}(\{\pi(p)\}\times\{f^1_\beta(p)\}\times...\times W_j\times...\times \{f^k_\beta\})$$
with adapted analog properties. 
\\

\footnotesize{Remark : A simple idea to define a multi-fiber structure based on the splitting of the fibers could have been to consider on the initial $\F$-fibration $\pi : \M\rightarrow \B\M$, with $\F=S\times W$, an additional map $\Phi=(h,f) :\M\rightarrow S\times W$ such that, for any $p\in\M$, the restriction 
$\Phi|_{\F_{\B p}} : \F_{\B p} \stackrel{\simeq}{ \rightarrow} S\times W$
 of $\Phi$ to the $\pi$-fiber $\F_{\B p}:=\pi^{-1}(\pi(p))$ is a diffeomorphism ;  hence the following diagram :
 \begin{align}
& \M  \stackrel{\Phi}{\longrightarrow}S \times W \nonumber \\
\pi & \downarrow \nonumber \\
& \B\M \nonumber
\end{align}
Then, simply define $ S_p := (f|_{\F_{\B p}})^{-1}(f(p))$ and $ W_p := (h|_{\F_{\B p}})^{-1}(h(p))$. But it is easy to see that in fact this gives a trivial fibration in the sense that $\M$ is then diffeomorphic to $\B\M \times S \times W$ : just consider $\phi : \M \rightarrow \B\M \times S \times W$, $x \mapsto (\pi (x), (h(x),f(x)))=(\pi(x),\Phi(x))$ whose inverse is 
$\phi^{-1} : (a,b)\mapsto (\Phi|_{\pi^{-1}(a)})^{-1}(b)$.
}
\normalsize

\subsubsection{Adapted charts.} 
As they are useful  to understand the situation, let us see how we get adapted charts, and what they look like. We take here $S=S^1$ as we will be mostly interested in this case. Consider $\M$ with a $(S^1\times W)$-multi fiber structure. As we saw above, starting with a trivialization chart of the fibration $\pi$ and composing with the splitting of the fiber  $ \F_{p}$ by $\psi_p$, we have in the neighborhood of any fixed point $p\in\M$ an adapted diffeomorphism of the following form:
\begin{align}
 \U_p  \stackrel{\phi_\U}{\longrightarrow} \, \, & \B\U_p \times (S^1\times W)_p \nonumber \\
& Id \downarrow \,\,\, \, \quad \quad \,\, \downarrow \psi_p  \nonumber \\
& \B\U_p \times \, S^1_p\times \, W_p \nonumber
\end{align}
As $S^1_p$ and $W_p$ are diffeomorphic to $S^1$ and $W$ respectively, we can now take coordinates $(x^i)$ on $\B\U_p$, $(u)$ on some neighborhood $\dot S^1_p$ of $p$ in $S^1_p$, and $(w^k)$ on some neighborhood $\W_p$ of $p$ in $W_p$, to obtain a chart of the form :
\begin{align}
 \U_p  \stackrel{\phi_\U}{\longrightarrow} \, \, & \B\U_p \times (S^1\times W)_p \nonumber \\
& Id \downarrow \,\,\, \, \quad \quad \,\, \downarrow \psi_p \nonumber \\
& \B\U_p \times \, \dot S^1_p\times \, \W_p \nonumber \\
& \downarrow \quad \quad \downarrow \quad \quad \downarrow   \nonumber \\
& (x^i, \quad  \, u, \,\quad \, w^k) \nonumber
\end{align}
Centering the chart so that the coordinates of $p$ are $(0,...,0)$, the coordinates expression of $f|_{\F_{p}}$ and $h|_{\F_{p}}$ are :
$$f|_{\F_{p}}(0,...,0, u,w^1,...,w^m)=(w^1,...,w^m)$$
$$h|_{\F_{p}}(0,...,0, u,w^1,...,w^m)=(u)$$
Indeed, they are submersions ! It  is then clear that :
$$(f|_{\F_{p}})^{-1}(f(x))=\{(0,..0,u,0,...,0)\}:=\dot S^1_p$$
$$(h|_{\F_{p}})^{-1}(h(x))=\{(0,...,0,0,w^1,...,w^m)\}:=\W_p$$
for some neighborhoods $\dot{S^1_p}$ and $\W_p$ of $p$ in $S^1_p$ and $W_p$ respectively, where the coordinates are defined. 

\subsubsection{Multi-fiber manifolds.} 
Our definition of the multi-fiber bundle structure is based on the fiber bundle atlas. We can also define directly on a manifold a mean to get well-defined fibers at each point, without referring to an existing bundle structure. This was the original idea of the second author, Michel Vaugon.
\begin{greybox}
Let $\M$ be a differential $n$-dimensional manifold, and let $S$ and $W$ be two compact manifolds of respective dimension $k$ and $l$.
A diffeomorphism : $$\phi:\U \rightarrow \B\U\times S \times W$$ where $\U$ is an open set in $\M$ and $\B\U$ an open set in  $\R^{n-k-l}$, will be called an \defin{observation diffeomorphism}, and the couple $(\U,\phi)$ an \defin{observation chart} (if $\B\U$ is understood). We note $\phi=(f^1,f^2,f^3)$ the threee components of $\phi$.

We say that $\M$ is a \defin{multi-fiber manifold} with fibers $S$ and $W$ if there exists a \defin{$(S,W)$-observation atlas}, that is, a family $\{(\U_\alpha, \phi_\alpha)\}_{\alpha \in A}$ of observation charts such that $\cup\U_\alpha=\M$, and satisfying for any $\alpha, \beta \in A$ and any $p\in \U_\alpha \cap \U_\beta$  :
\begin{itemize}
	\item $\phi_\alpha^{-1}(\{ f^1_\alpha(p)\}\times S\times \{f^3_\alpha(p)\})=\phi_\beta^{-1}(\{f^1_\beta(p)\}\times S\times \{f^3_\beta(p)\})$
	\item $\phi_\alpha^{-1}(\{f^1_\alpha(p)\}\times\{f^2_\alpha (p)\}\times W)=\phi_\beta^{-1}(\{f^1_\beta(p)\}\times\{f^2_\beta(p)\}\times W)$
\end{itemize}
The definitions, for any $p\in \M$, of fibers $S_p$ and $W_p$ is then given as above.	
\end{greybox}
An observation atlas can be completed in a \emph{complete observation atlas} in the same manner as for a classical differential manifold atlas.

Whereas the multi-fiber bundle structure is a generalization of the fiber bundle structure (or fibration), the multi-fiber manifold structure can be seen as a generalization of the \emph{foliation} structure.

\emph{Remark and proposition : if we consider a single compact fiber W, and an observation atlas of observation charts of the form $\phi:\U \rightarrow \B\U\times W$ satisfying the adapted compatibility condition : $\phi_\alpha^{-1}(\{f^1_\alpha(p)\}\times W)=\phi_\beta^{-1}(\{f^1_\beta(p)\}\times W)$,
it is fairly easy to prove that $\M$ can be equipped with a fiber bundle structure $\pi : \M\rightarrow \B\M$ with fiber $W$ for some manifold $\B\M$. Indeed, consider $\sim$ defined by $p\sim p'$ if $\phi_\alpha^{-1}(\{f^1_\alpha(p)\}\times W)=\phi_\alpha^{-1}(\{f^1_\alpha(p')\}\times W)$, that is, if $W_p=W_{p'}$. Define $\B\M:=\M/\sim$. The compacity of $W$ ensures $\B\M$ is Hausdorf ; if $W$ is not compact, one must add the requirement that $\phi_\alpha^{-1}(\{f^1_\alpha(p)\}\times W)$ is closed in $\M$ for all $\alpha$ and all $p$.}
\\

The main difference between this structure and the multi-fiber bundle structure, is that it makes no reference to a "natural" manifold $\B\M$ linked to the "horizontal" distribution $H_p=(S_p \times W_p)^\bot$. The use of the bundle structure can bring formulas \emph{à la O'Neill} linking the geometries of $\B\M$, $S$ and $W$. On a heuristic point of view, keeping the bundle structure keep the idea of "small compact dimensions" attached to classical 4-dimensional spacetime $\B\M$, whereas choosing as model a multi-fiber manifold with a $(S,W)$-observation atlas is  more radical as it makes no reference to a specific 4-dimensional manifold.

\subsubsection{Construction.}
We imitate here the classical constructions of fiber bundles using cocyles with value in a subgroup of the diffeomorphisms group of the fiber. This starts by observing that the overlap maps $\phi_\alpha \circ\phi_\beta^{-1} : (\B\U_\alpha \cap \B\U_\beta)\times\F \rightarrow (\B\U_\alpha \cap \B\U_\beta)\times\F$ give rise to diffeomorphisms of the fiber $\F$ : $$\Phi_{\alpha\beta}(p):=\Phi_\alpha |_{\F_p}\circ \Phi_\beta |_{\F_p}^{-1} :\F\rightarrow \F, $$ where as always, $\phi$ is written $\phi=(\pi,\Phi)$.
In our multi-fiber case, with $\F=S\times W$, we have two more diffeomorphisms, writing again $\Phi_\alpha=(h_\alpha, f_\alpha)$ :
$$h_{\alpha\beta}(p):=h_\alpha |_{S_p}\circ h_\beta |_{S_p}^{-1} :S\rightarrow S $$
$$f_{\alpha\beta}(p):=f_\alpha |_{W_p}\circ f_\beta |_{W_p}^{-1} :W\rightarrow W. $$
Indeed, $h_\alpha |_{S_p}:S_p\rightarrow S$ and $f_\alpha |_{W_p}:W_p\rightarrow W$ are diffeormorphisms. We therefore have, for each $\alpha,\beta$ such that $\B\U_\alpha \cap \B\U_\beta\neq \emptyset$,  maps $p\mapsto h_{\alpha\beta}(p) \in Diff(S)$ and $p\mapsto f_{\alpha\beta}(p) \in Diff(W)$.
These diffeomorphisms satisfy a cocycle relation, in the sense that : 
$$h_{\alpha\alpha}(p)=Id_S\,\,,\,\, h_{\alpha\beta}(p)=h_{\beta\alpha}(p)^{-1}\,\,,\,\, h_{\alpha\beta}(p)\circ h_{\beta\gamma}(p)=h_{\alpha\gamma}(p),$$ and similarly for $f_{\alpha\beta}(p)$. 
These cocycles are the building blocks for our multi-fiber structure.\\

So let $\B\M, S,W$ be three manifolds, $S$ and $W$ being compact. Let $H$ and $F$ be two Lie groups that acts on the left on $S$ and $W$ respectively. Let be given an open cover $(\B\U_\alpha)_{\alpha\in A}$ of $\B\M$. A $H$-cocycle for $(\B\U_\alpha)$ is the assignment of a smooth map $h_{\alpha\beta}:(\B\U_\alpha \cap \B\U_\beta)\rightarrow H$ to every nonempty intersection $\B\U_\alpha \cap \B\U_\beta$ such that the cocycle conditions holds for $h_{\alpha\beta}$. We suppose we are given an $H$-cocycle $(h_{\alpha\beta})$ and a $F$-cocycle $(f_{\alpha\beta})$, and from now on, to make the writings simpler, we do as if $H$ and $F$ were subgroups of the diffeomorphisms groups $Diff(S)$ and $Diff(W)$ respectively, with the natural actions ; the reader will easily adapt what follows to the more general case of Lie groups acting on $Diff$.

With these datas, we can construct our $(S,W)$-multi-fiber strucure over $\B\M$. For this, we consider the disjoint union
$$\Sigma:=\coprod_\alpha\, \{\alpha\}\times\B\U_\alpha\times S \times W$$
and on $\Sigma$, the equivalence relation $\sim$ defined by $$(\alpha, p, x,y)\sim (\beta, p',x',y')\quad \mbox{iff}$$ $$p=p',\quad x'=h_{\alpha\beta}.x,\quad \mbox{and} \quad y'=f_{\alpha\beta}.y$$
We define $\M:=\B\M/\sim$ and $\pi :\M\rightarrow\B\M, \, [\alpha, p,x,y]\mapsto p$, where $ [\alpha, p,x,y]$ is the equivalence class of $(\alpha, p,x,y)$. 
Now, defining $\Phi_{\alpha\beta}:=(h_{\alpha\beta},f_{\alpha\beta}):\F\rightarrow\F$ for $\F:S\times W$, we obtain a $H\times F$-cocycle acting on $\F$. Thanks to this, it is classical to prove that $\M$ has a $\F$-bundle structure.

Indeed, an bundle atlas on $\M$ is obtained by defining $\U_\alpha=\pi^{-1}(\B\U_\alpha)$ and $\phi_\alpha :\U_\alpha\rightarrow\B\U_\alpha\times S \times W$ as the map $m\in \U_\alpha\mapsto(p,x,y)$ such that $(\alpha,p,x,y)\in m$ ; remember that $m$ is an equivalence class. This is well defined because, as $h_{\alpha\alpha}=Id_S$ and $f_{\alpha\alpha}=Id_W$, $(\alpha,p,x,y)\sim(\alpha,p',x',y')$ iff $p=p'$, $x=x'$ and $y=y'$. The family $(\B\U_\alpha, \phi_\alpha)_{\alpha\in A}$ then constitutes a $\F$-bundle atlas.

We won't go in the detailed proof of this. We just indicate what an overlap map $\phi_\alpha \circ \phi_\beta^{-1}$
looks like. So let $\B\U_\alpha\cap\B\U_\beta\neq\emptyset$, and $p\in \B\U_\alpha\cap\B\U_\beta$. $\phi_\beta^{-1}(p,x,y)=[\beta,p,x,y]\in\pi^{-1}(\B\U_\alpha)$ thus $[\beta,p,x,y]=[\alpha,q,x',y']$ which means $p=q$, $x'=h_{\alpha\beta}.x$ and $y'=f_{\alpha\beta}.y$ ; i.e. $\phi_\beta^{-1}(p,x,y)=[\alpha,p, h_{\alpha\beta}.x, f_{\alpha\beta}.y]$. Therefore $\phi_\alpha \circ \phi_\beta^{-1}(p,x,y)=(p, h_{\alpha\beta}.x, f_{\alpha\beta}.y)=(p,\Phi_{\alpha\beta}(x,y))$. Setting a smooth structure on $\M$ and proving that it is Hausdorf is classic.
\\

Remark : A general diffeomorphism $\Phi=(h,f) : S\times W \rightarrow S\times W$ is written $\Phi(x,y)=(h(x,y),f(x,y))$. Here, the cocycle $\Phi_{\alpha\beta}$ are of the special form : $\Phi_{\alpha\beta}(p)(x,y)= (h_{\alpha\beta}(p)(x),f_{\alpha\beta}(p)(y))$. \\

Now, we need to check that we can extract from the $\F$-bundle atlas $(\B\U_\alpha, \phi_\alpha)_{\alpha\in A}$ a $(S,W)$-multi-fibers atlas. Once again, we write, $\phi_\alpha=(\pi,\Phi_\alpha)=(\pi,h_\alpha,f_\alpha)$.
  We have to prove that for any $m\in \pi^{-1}(\B\U_\alpha \cap \B\U_\beta)$ we have
$\phi_\alpha^{-1}(\{\pi(m)\}\times S\times \{f_\alpha(m)\})=\phi_\beta^{-1}(\{\pi(m)\}\times S\times \{f_\beta(m)\})$ and similarly for $W$. 

If $m\in \pi^{-1}(\B\U_\alpha \cap \B\U_\beta)$, it can be written $m=[\alpha, p,x,y]=[\beta, p,x',y']$ with $p=\pi(m)$, $x'=h_{\alpha\beta}.x$, $y'=f_{\alpha\beta}.y$. If $u\in \phi_\alpha^{-1}(\{\pi(m)\}\times S\times \{f_\alpha(m)\})$, $u$ is written $u=[\alpha,p',a,b]$. But $\pi(u)=\pi(m)$ implies $p=p'$ by definition of $\pi$. Similarly, by definition of $\phi_\alpha$, $\phi_\alpha(u)=(p,a,b)$ and $\phi_\alpha(m)=(p,x,y)$, so $f_\alpha(u)=f_\alpha(m)$ implies $b=y$. Therefore $u=[\alpha, p, a, y]$. Now, as $\pi(u)=\pi(m)$, $u\in \pi^{-1}(\B\U_\alpha \cap \B\U_\beta)$. So $u=[\beta, p, a',b']$ with $a'=h_{\alpha\beta}.a$ and $b'=f_{\alpha\beta}.y$. By definition of $\phi_\beta$, $\phi_\beta (u)=(p,a',b')$. But $a'=h_{\alpha\beta}.a\in S$ and $b'=f_{\alpha\beta}.y=y'=f_\beta(m)$. Therefore, $\phi_\beta(u)\in \{\pi(m)\}\times S\times \{f_\beta(m)\}$.

We just proved that $\phi_\beta^{-1}(\{\pi(m)\}\times S\times \{f_\beta(m)\})\subset \phi_\alpha^{-1}(\{\pi(m)\}\times S\times \{f_\alpha(m)\})$. Exchanging the role of $\alpha$ and $\beta$ above, we get that $\phi_\alpha^{-1}(\{\pi(m)\}\times S\times \{f_\alpha(m)\})=\phi_\beta^{-1}(\{\pi(m)\}\times S\times \{f_\beta(m)\})$. We proceed analogously to prove that $\phi_\alpha^{-1}(\{\pi(m)\}\times\{h_\alpha(m)\}\times W)=\phi_\beta^{-1}(\{\pi(m)\}\times\{h_\beta(m)\}\times W)$.

Therefore, $(\B\U_\alpha, \phi_\alpha)_{\alpha\in A}$ is also a $(S,W)$-multi-fiber bundle atlas. (Note however that $(\B\U_\alpha, \phi_\alpha)_{\alpha\in A}$ might be completed in a larger $ç(S\times W)$-bundle atlas).

\subsubsection{Building objects on fibers.}
Conversly, let $(\B\U_\alpha, \phi_\alpha)_{\alpha\in A}$ be a $(S,W)$-multi-fiber bundle atlas for a manifold $\M$, coming from a $\F$-bundle atlas for the fibration $\pi : \M \rightarrow\B\M$ where $\F=S\times W$. For each $p\in \M$, the fibers $S_p$ and $W_p$ are well defined. Again, we write $\phi_\alpha=(\pi,\Phi_\alpha)=(\pi,h_\alpha,f_\alpha)\,$. Then, 
for any $\alpha,\beta$ such that $\pi(p)\in \B\U_\alpha\cap\B\U_\beta$, $h_\alpha|_{S_p}$ and $h_\beta|_{S_p}$ are diffeomorphisms $S_p\rightarrow S$, and $f_\alpha|_{S_p}$ and $f_\beta|_{S_p}$ are diffeomorphisms $W_p\rightarrow W$.

For any $\alpha,\beta$ with $\B\U_\alpha\cap\B\U_\beta\neq \emptyset$ we therefore have maps :
$$h_{\alpha\beta} : p\mapsto h_{\alpha\beta}(p):=h_\alpha|_{S_p}\circ h_\beta|_{S_p}^{-1}\in Diff(S)$$
and
$$f_{\alpha\beta} : p\mapsto f_{\alpha\beta}(p):=f_\alpha|_{W_p}\circ f_\beta|_{W_p}^{-1}\in Diff(W).$$
For each $p\in\M$, $H_p:=(h_{\alpha\beta}(p))$ and $F_p:=(f_{\alpha\beta}(p))$ define cocycles in $Diff(S)$ and $Diff(W)$ respectively. We consider the subgroups $H$ of $Diff(S)$ and $F$ of $Diff(W)$ generated by the $H_p$ and $F_p$ respectively : $H:=<(\cup_{p\in\M} H_p>)$ and $F:=<(\cup_{p\in\M}F_p)>$.

We suggest here some possible ways to build tensors on $\M$ out of similar objects defined on $S$ or $W$, that could not be defined using the sole $\F$-bundle structure of $\M$. (We will use the fiber $S$, but obviously the same constructions can be made using $W$). \\

\defin{Functions :} Let $f\in\C^\infty(S)$. Let's say $f$ is $H$-invariant on $S$ if $\forall h\in H$, $h^*f=f$, that is, $f\circ h=f$. We can then define a smooth function $\tilde f$ on $\M$ by setting $\tilde f(p)=f(h_\alpha(p))$ for any multi-fiber chart $(\B\U_\alpha, \phi_\alpha)$, such that $p\in\U_\alpha:=\pi^{-1}(\B\U_\alpha)$.

The idea if of course that if $p\in\U_\beta$, because $h_\alpha|_{S_p}\circ h_\beta|_{S_p}^{-1}=h_{\alpha\beta}(p)\in H$, we have $f(h_\beta(p))=f((h_\alpha|_{S_p}\circ h_\beta|_{S_p}^{-1})(h_\beta(p))=f(h_\alpha(p))$.

We could not build such a function on $\M$ from a function on $S$ with the sole $\F$-bundle structure.
\\

\defin{Vector fields :} Let $X\in \Gamma (TS)$ be a vector field on $S$. Let's say $X$ is $H$-invariant on $S$ if $\forall h\in H$ $h_* X=X$, that is, $\forall x\in S$, $T_x h(X(x))=X(h(x))$. 

Let $p\in\M$ with $p\in\U_\alpha$ for some $(\B\U_\alpha, \phi_\alpha)$ of the multi-fiber atlas. Then $h_\alpha|_{S_p}:S_p\rightarrow S$ is a diffeomorphism. We define $\tilde X(p)= (T_p h_\alpha|_{S_p})^{-1}(X(h_\alpha(p))$.

If $p\in\U_\beta$, we have, as $h_\alpha|_{S_p}\circ h_\beta|_{S_p}^{-1}=h_{\alpha\beta}(p)\in H$ :
\begin{align}
	(T_p h_\beta|_{S_p})^{-1}(X(h_\beta(p))&=(T_p h_\beta|_{S_p})^{-1}(X(h_\beta|_{S_p}\circ h_\alpha|_{S_p}^{-1}(h_\alpha(p))\\
	&=(T_p h_\beta|_{S_p})^{-1}T_{h_\alpha(p)}(h_\beta|_{S_p}\circ h_\alpha|_{S_p}^{-1})(X(h_\alpha(p))\\
	&=(T_p h_\alpha|_{S_p})^{-1}(X(h_\alpha(p))
\end{align}
So $\tilde X$ is a well-defined vector field on $\M$ that, here again, could not have been defined from a vector field $X$ on $S$ without the multi-fiber structure.
\\

\defin{The case of covariant tensors :} Pulling back (covariant) tensors on $\M$ from tensors on $S$ is more delicate, as $h_\alpha$ only induces an isomorphism on the subspace $T_pS_p$ of the whole tangent space $T_p\M$. So even though we can push forward a vector $v$ of $T_p\M$ using $T_p h_\alpha$, $T_{h_\alpha(p)}(h_\beta|_{S_p}\circ h_\alpha|_{S_p}^{-1})(T_p h_\alpha(v))$ will not be well-defined outside $T_pS_p$.

It appears that the only way to define properly a tensor $\tilde A$ on $\M$ from the tensor $A$ on $S$ is to suppose that we are given an \emph{horizontal distribution} $\Hor$ on $\M$, that is, a smooth family $(H_p)_{p\in\M}$ of subspaces of $T_p\M$, such that at each $p\in\M$, $T_p\M=H_p\oplus T_p\F$ ; in this case , we also have $\T_p\M=H_p\oplus T_pS\oplus T_pW$.

For example, if $\M$ comes equipped with a Riemannian metric $\g$, an obvious choice is to take $H_p:=( T_pS\oplus T_pW)^\bot$. If $\g$ is semi-Riemannian, we have to require appropriate signature compatibility on the fibers to ensure that $H_p$ so defined is a supplementary to $ T_pS\oplus T_pW$ in $\T_p\M$.

We therefore suppose now that we are given such an horizontal distribution $\Hor$. Any tangent vector $v\in T_p\M$ can be written uniquely $v=v_H+v_S+v_W$, with $v_H\in H_p$, $v_S\in T_pS_p$ and $v_W\in T_pW_p$.

So let's $A\in T^{2,0}(S)$ be a covariant 2-tensor field on $S$. 
We shall say that $A$ is $H$-invariant on $S$ if, $\forall h\in H$, $h^*A=A$, that is : $A_{h(x)}(T_xh.u,T_xh.v)=A_x(u,v)$, $\forall x\in S$, $\forall u,v \in T_xS$.

Let again $p\in\M$ with $p\in\U_\alpha$ for some $(\B\U_\alpha, \phi_\alpha)$ in the multi-fiber atlas, and let $u,v\in T_p\M$. We define a (2,0)-tensor field $\tilde A$ on $\M$ by setting :
$$\tilde A_p (u,v)=((h_\alpha|_{S_p})^*A)_p(u_S,v_S):=A_{h_\alpha(p)}(T_p h_\alpha|_{S_p})(u_S),T_p h_\alpha|_{S_p})(v_S))$$
We check, as above in the case of a vector field, that $\tilde A$ is a well-defined (2,0)-tensor field on $\M$.

 \subsubsection{Building Metrics.}
 A natural question when given a classical fiber bundle $\pi:\M\rightarrow B$ over a manifold base $B$ and fiber $\F$, and metrics $\g_B$ and $\g_\F$ on $B$ and $\F$ respectively, is to build a natural metric $\g$ on $\M$ out of $\g_B$ and $\g_\F$.

However, just as we saw above for any covariant tensors, the difficulty to pull back tensors from $B$ or $\F$ is due to the absence of a canonical supplementary space to $T_p \F$ in $T_p\M$. 

As this natural question extends naturally to our multi-fiber structure case, and as a solution lies on the same requirement (which is the existence of a given horizontal distribution), we address it here in this more general case.

So let $\pi : \M\rightarrow\B\M$ be a $(S,W)$-multi-fiber bundle, and let be given the same cocycles and diffeomorphisms groups datas as at the beginning of the previous section.

We suppose that $\B\M$ is equipped with a metric $\g_B$ ($B$ for Base !), that $S$ is equipped with a $H$-invariant metric $\g_S$, and $W$ with a $F$-invariant metric $\g_W$. (Invariance here is understood as in the above case of covariant tensors in the previous section).
Define the vertical space $\V_p:=Ker (T_p\pi)$, let $\V=\cup_p \V_p$ be the \emph{vertical distribution}. We suppose we are given, at each $p\in\M$, a supplementary space $H_p$ to $\V_p$, called \emph{horizontal space}, and that the distribution $\Hor:=\cup_p H_p$ is a smooth distribution. At each $p\in\M$, we have :
$$T_p\M=\V_p\oplus H_p$$
$$\V_p=T_p\F_p=T_pS_p\oplus T_pW_p$$
And therefore : $T_p\M=H_P\oplus T_pS_p\oplus T_pW_p$.

Any $v\in T_p\M$ can thus be written : $v=v_H+v_S+v_W$, with $v_H\in H_p$, $v_S\in T_pS_p$ and $v_W\in T_pW_p$.

It is now easy to define a metric $\g$ on $\M$. Letting $p$ be any point in $\M$ with $p\in\U_\alpha$ for some $(\B\U_\alpha, \phi_\alpha)$ in the multi-fiber atlas, and $u,v\in\T_p\M$, define :
$$\g_p(u,v):=\pi^*\g_B(u_H,v_H)+(h_\alpha|_{S_p})^*\g_S(u_S,v_S)+(f_\alpha|_{W_p})^*\g_W(u_W,v_W).$$

Remember our notation from the previous section : if $\w t$ is a covariant tensor on $S$ or $W$, $\w {\tilde t}$ is the tensor pulled back by $h_\alpha|_{S_p}$ or $f_\alpha|_{W_p}$ respectively. We can also write the above definition of $\g$ :
$$\g:=\pi^*\g_B+\w{\tilde \g}_S+\w{\tilde \g}_W.$$
It is also now very easy to imagine \emph{warped metrics} in the spirit of O'Neill warped product : letting $a,b\in\C^\infty(\B\M)$ be positive functions, set :
$$\g:=\pi^*\g_B+(a\circ\pi)^2.\w{\tilde \g}_S+(b\circ\pi)^2.\w{\tilde \g}_W.$$

O'Neill type formulae and results can then be obtained linking geodesics or Ricci curvature (for example) of $(\M,\g)$ to the geometries of $(\B\M,\g_B)$, $(S,\g_S)$ and $(W,\g_W)$.

More generally, in [19], Michel Vaugon uses a metric conformal to the above, where the conformal factor function, $f\in\C^\infty(\M)$, $f>0$, is used to model quantum phenomena :
$$\g:=f^2.\pi^*\g_B+(a\circ\pi)^2.\w{\tilde \g}_S+(b\circ\pi)^2.\w{\tilde \g}_W.$$

\subsubsection{The special case of $S^1$, electromagnetic potential, and an example :} 
For our $(S,W)$-multi-fiber bundle structure, the case $S=S^1$, $S^1$ being the standard circle,  will be particularly important. Indeed, it is the fiber to be used to include electromagnetism in the geometric frame of general relativity according to Kaluza-Klien idea.

Besides, it is a very tractable case, as in particular we can define easily on $\M$ a natural vector field $\w Y$ associated to each fiber $S_p$ without using all the machinery of cocycles and invariance ; $\w Y$ will of course be the electromagnetic potential.

Indeed, consider $\M$, a $(S^1,W)$-multi fiber bundle. A circle fiber $S^1_p$ is well defined at each $p\in\M$. We suppose the multi-fiber structure on $\M$ is compatible with the natural orientation of $S^1$ (considering for example $S^1$ as $\R / \Z$). Then we can define unambiguously $\w Y$ to be the vector field defined at each $x \in \M$ to be tangent to the fiber $S^1_x$ and such that $\g (\w Y,\w Y)=-1$, with the chosen orientation for $S^1$. \\

As we already said, the main idea leading to the multi-fiber structure is that we will use the other fibers $W_p$ to model other physical interactions, but still keeping the possibility to use the geometry of the total fiber $\F_p \simeq S^1_p \times W$.

As an example we cite the use by the second author, Michel Vaugon, of a 
\defin{$(S^1\times S^3)$-multi-fiber structure on a manifold $\M$.}

The 3-dimensional sphere $S^3$ is the classical geometric space used in quantum theory to describe \emph{spin}. Indeed, $S^3$ carries a natural frame field as well as canonical endomorphisms giving precisely the spin matrixes. This lead Michel to use a multi-fiber bundle structure on a manifold $\M$ with fiber $\F=S^1\times S^3$ to give a geometric model of electromagnetic \emph{and} spin effects on $\M$. The $S^1$ component was used to define the electromagnetic potential, and the $S^3$ component to define objects related to the spin. Furthermore, the modelization of the physical effects of both electromagnetism and spin, for instance to describe the Stern and Gerlach experiment, required to use the geometry of the \emph{total} fiber  $\F=S^1\times S^3$, so the full multi-fiber structure was used. 

More precisely, this $S^1\times S^3$-multi-fibers structure, and the objects described below that can be built with it,  was the geometric starting point for the second author to build a new approach towards a geometric unification of General Relativity and Quantum Physics. It is based on special metrics on a $(5+k)$-dimensional manifold modelizing quantum particles physics in spacetime. In this setting, only the metric is relevant, no objects or laws are added, these appear as geometric quantities issued from curvature and geometric theorems, such as Bianchi identity, linked with the objects described below. See [19].

We therefore consider that the spacetime $(\M,\g)$ is a $(S^1,S^3)$-multi-fibers bundle.
At each point $x\in\M$, we then have two naturally defined fibers : $S^1_x$ and $S^3_x$.

As before, $S^1$ gives the electromagnetic potential $\w Y$ :  $\w Y$ is the vector field defined at each $x \in \M$ to be tangent to the fiber $S^1_x$ and such that $g (Y,Y)=-1$, with the chosen orientation for $S^1$.

Then, on $S^3$, any function, vector field, or covariant tensor, invariant by a subgroup of $Diff(S^3)$ containing the cocycles induced by the overlap maps of the multi-bundle atlas can be used to construct analog objects on $\M$. 
In his work for example, Michel uses the spectral theory of $S^3$ to transpose a Hilbertian basis of $L^2$ functions from $S^3$ to $\M$.

These objects could not be defined with a simple $(S^1\times S^3)$-fiber bundle structure on $\M$.

\section{Extending Kaluza-Klein model of spacetime.}
\subsection{Multi-fiber Kaluza-Klein spacetime.}
We propose as example a possible model for space-time. For illustration of the possibilities, and further developments by Michel Vaugon to be found in [ref : New approach to Kaluza-Klein theory, and A Mathematicians’ View of Geometrical Unification of General Relativity and Quantum Physics], we choose the fifth dimension to be timelike. The reader uncomfortable with the two timelike dimensions, can still consider the signature on $S^1$ to be spacelike, only minor sign changes will be required in front of expressions using $e$ or $Y$, but all what follows remains in fact essentially unchanged.
\begin{greybox} 
Let $S^1$ be the classical circle with a chosen orientation, and let $W$ be a compact manifold of dimension $m$. 

A \emph{Spacetime} is a semi-Riemannian manifold $(\M,\g)$ of dimension $5+m$ equipped with a multi-fiber bundle structure, of fiber $\F=S^1\times W$ as defined above. 

We suppose that the metric  $\g$ is compatible with the multi-fiber bundle structure, $\g$ being of  signature $(-1)$ on $S^1$, and $(+,...,+)$ on $W$ ; the total signature of $\g$ is thus $(-,+,+,+,-,+,...,+)$. $\B\M$ is "classical" spacetime.

The \defin{horizontal space} at a point $x$ is $H_x := (\T_x(S^1\times W)_{\B x})^\bot $, and $\g_x$ is of signature $(-,+,+,+)$ on $H_x$. $H_x$ represents the local and classical Minkowski spacetime at $x$.
\end{greybox}

The effect of the "extra" $m$ dimensions carried by $W$ will be modeled via the geometry of the map $\pi: \M \rightarrow \B \M$ and $(S^1, W)$-multi-fiber structure,  and via the metric $\g$ or its Einstein curvature $\w G$. The idea is also that what passes to the quotient can be neglected.
\\

We can also add to the definition of $\g$ being compatible with the given $(S^1,W)$-multi-fiber structure on $M$ the following requirement : for any pair of adapted charts $\phi_i$, $\phi_j$ as defined above, $\forall x \in \U_i\cap \U_j$, 
$\phi_i^*(\partial_t)_x$ and $\phi_j^*(\partial_t)_x$ are timelike and in the same time orientation, i.e. $g(\phi_i^*(\partial_t)_x,\phi_j^*(\partial_t)_x)<0$, where $\partial_t$ is the tangent vector to the canonical coordinates $(t,x,y,z)$ on $\Theta_i \subset \mathbb{R}^4$. This condition gives a "classical" time-orientation on every apparent space-time $H_x$, varying differentially with $x$.
\\

\textbf{The Electromagnetic potential.}

We suppose that we have the canonical standard orientation on $S^1$. We can then define:
\begin{greybox}
We define $\w Y$ to be the vector field defined at each $x \in \M$ to be tangent to the fiber $S^1_x$ and such that $\g  (\w Y,\w Y)=-1$, with the chosen orientation for $S^1$ ; $\w Y$ is called the \defin{electromagnetic potential}. We then define the differential  $\w F:=\dd (\w Y^\flat)$ where 
$\w Y^\flat$ is the 1-form associated to $\w Y$ by $\g$ ; $\w F$ is the \defin{electromagnetic field}. We always suppose from now on that $\w Y$ is a Killing vector field. (Note that if $\w Y$ is Killing and of constant norm, it is necessarily geodesic.) The local diffeomorphisms generated by $\w Y$ are therefore isometries.
\end{greybox}

The next proposition, easy to prove, shows why it is natural to suppose that $\w Y$ is a Killing vector field, when supposing that the compact dimensions are "small":
\begin{greybox}
\begin{proposition}
\textbf{Averaging the metric on $S^1$ :} Let $\w Y$ be tangent to the fiber $S^1_x$ and such that $g (Y,Y)=-1$ as above, but without supposing that $\w Y$ is a Killing vector field. Let $\sigma$ be the 1-parameter group of diffeomorphisms associated to the flow of $\w Y$. Define the "averaged" metric $\B \g$ by :
$$\forall x\in\M , \quad \B\g_x:=\frac{1}{\ell_x}\int_{t_0}^{t_0+\ell_x}(\sigma^*(t)\g)_x.dt $$
where $\ell_x$ is the length of $S^1_x$ relative to $\g$. ($\B\g_x$ does not depend on the choice of $t_0$ as $\sigma_x(.)$ is periodic, of period $\ell_x$.)  Then, $\B\g (\w Y,\w Y)=-1$, and $\forall s\in \R$, $\sigma^*(s).\B\g=\B\g$. That is, $\w Y$ is a Killing vector field for $\B\g$.
\end{proposition}
\end{greybox}

\subsection{Application to an electrically charged fluid}

To illustrate possible use of our multi-fibers structure, we give an application to a new approach to Kaluza-Klein theory. Details and proofs can be found in [3]. We start by considering the inclusion of electromagnetism in the classical frame of General Relativity. Then we show how the multi-fibers structure gives a possible way to model or encode deviations from standard 4-dimensional General Relativity, or "dark" effects such as dark matter or energy.

We note $\w G=\w{Ric}-\frac{1}{2}\w S.\g$ the Einstein curvature. $^e\w G$ is the associated endomorphisms field. We then define $^e \w G_H$, the endomorphisms field on the horizontal subspaces $H_x$, defined by $^e \w G_H=pr_H \circ (^e \w G_{|H})$, where for $x \in \M$, $(pr_H)_{|x}$ is the orthogonal projection of $\T_x \M$ on $H_x$. This tensor will be very important to define fluids in our extended Kaluza-Klein spacetime.
\\

\textbf{Dust charged matter fluid }
\begin{greybox}
\begin{definition}
 A domain $\Omega \subset \M$ is a \defin{dust charged matter fluid domain} if and only if its Einstein curvature tensor can be written  at every point 
 $$\w G=\mu X \otimes X + \alpha \w Y \otimes \w Y$$
with the condition that, at each point $x$, $pr_H(X)$ is a basis for a timelike 1-dimensional eigenspace of $^e\w G_H$ of eigenvalue $-\mu<0$. $X$ is then unique for this decomposition.

\textbf{Associated "classical" data : } For such a perfect fluid without pressure, there is a unique decomposition (once a time orientation is chosen) :
$$\w G= \mu X_0 \otimes X_0 +e(X_0 \ts \w Y+\w Y \ts X_0)+\gamma \w Y \ts  \w Y,$$
where $\g(X_0,X_0)=-1$ and $X_0 \bot \w Y$. $\mu$ is called mass density, $e$ the charge density. These are canonically given by : $$\mu=-\w G(X_0,X_0)$$ $$e=\w G(\w Y,X_0)$$ $$\gamma= \w G(\w Y,\w Y)$$
 We then have 
 $$X=X_0+\frac{e}{\mu}\w Y$$
The vector field $X=X_0+ \frac{e}{\mu}\w Y$ is called the vector field of the fluid, and the associated flow, the flow of the fluid.
The vector field $X_0$ will be called the \emph{apparent, or visible,} field of the fluid, and the associated flow, the \emph{apparent, or visible,} flow. Note that at each point, by definition, $X_0(x)\in H_x$.
 \end{definition}
 \end{greybox}

The following theorem is a consequence of the purely geometric Bianchi identity applied to the Einstein curvature $\w G$. No "laws" need to be added.
\begin{greybox}
\begin{theorem}\textbf{Dynamics of charged dust.}
For the domain of a perfect charged fluid without pressure where $\w G$ and its associated classical data are written :
$$\w G=\mu X \otimes X + \alpha \w Y \otimes \w Y= \mu X_0 \otimes X_0 +e(X_0 \ts \w Y+\w Y \ts X_0)+\gamma \w Y \ts  \w Y,$$
Bianchi identity gives:
\begin{itemize}
\item Conservation Laws: $X_0(\frac{e}{\mu})=\Div (\mu X_0)=\Div (eX_0)=0$ 
\item Maxwell equations :  $\dd \w F=0$ and $(\Div \w F)^\sharp =2eX_0-(2\gamma +\w S_g)\w Y$
\item Free Fall : $X=X_0+\frac{e}{\mu}\w Y$ is a geodesic vector field.
\item Lorentz equation: $\mu \wnab_{X_0} X_0=e. ^e\w F(X_0)$. This is just free fall read on $H$.
\end{itemize}
When projected on the "classical" 4-dimensional space-time $H=\w Y^\bot$, these equations are the classical equations of physics. 

Note that Lorentz equation is obtained from the geodesic motion of $X=X_0+\frac{e}{\mu}\w Y$ by developing $\wnab_X X=0$ and writing (projecting) this equation on the horizontal space $H=\w Y^{\bot}$, noting that $\wnab_{X_0} X_0 \bot \w Y$, which means that $\wnab_{X_0} X_0\in H$. We therefore see that : 
\begin{center}
 \emph{Free fall for $X$ is equivalent to Lorentz equation for $X_0$}. 
 \end{center}
\end{theorem}
\end{greybox}
  (Remember that the first Maxwell equation, $\dd \w F=0$, is always obvious as we set $\w F=\dd \w Y^\flat$.)
 
 \begin{proof}
The proof is available in [3], but here is the scheme :

First some properties of $\w Y$ are established. Then, considering $\Div  \w G$ as a vector, i.e. identifying $\Div \w G$ and $(\Div \w G)^\sharp$, and noticing that $\Div \w G =0$ by Bianchi identity, one compute:

\begin{itemize}
\item $\g(\Div \w G, \w Y)$, this will be charge conservation law.
\item $\g(\Div \w G, X_0)$, this will be mass (baryonic number) conservation law.
\item $pr_H (\Div \w G)$, this be the equation of motion. In the case of dust $pr_H (\Div \w G)=\Div \w G$, so the equation of motion is simply $\Div \w G=0$.
\end{itemize}
Then, $\Div \w F=\Div (\dd\w Y^\flat)$ is computed, which gives the second Maxwell law, the first, $\dd\w F=0$, being obvious as $\w F=\dd (\w Y^\flat)$.

Finally, to prove that the flow of $X$ is geodesic, we simply compute $\wnab_X X$, noticing that $X_0 (\frac{e}{\mu})=0$ ; this leads to $\wnab_X X=0$.

 \end{proof}

It is very important to note that we have never mentioned any kind of energy-momentum tensor. The point is that from our point of view, this concept has no meaning. Indeed, here is our frame of ideas :

A/: Space-time is a five dimensional semi-Riemannian manifold satisfying definition 1

B/: Instead of defining an energy-momentum tensor, we caracterize a domain of space-time by a geometric type. We \emph{then} define physical concepts by geometric caracteristics of curvature.

C/:  Physical equations are projection of the Bianchi identity $\Div \w G=0$ on the 4-dimensional subspace $H=\w Y^{\bot}$ modelizing our classical 4-dimensional space-time.
\\

\textbf{General charged matter fluid }

In the result above, we only used an object, $\w Y$ defined on the $S^1$ component of the total fiber $\F=S^1\times W$. Now, with the same notations as for a dust fluid, we can define a general fluid to be a domain $D$ of $\M$ where the Einstein curvature can be written :
$$
\w G =\mu X \ts X+ \alpha \w Y\ts \w Y+P
        =\mu X_0 \ts X_0+ e(X_0 \ts \w Y+\w Y \ts X_0)+(\alpha+\frac{e^2}{\mu}). \w Y \ts  \w Y+P\, .$$
Here, the tensor $P$ is the pressure, or constraint, tensor of the fluid. If we note $\mathrm{T}_x$ the subspace of dimension 2  of $T_xM$ generated by $X_0(x)$ and $Y_x$, then $P$ is just  $P=\w G-\w G_{|\mathrm{T}_x}$, where $\w G_{|\mathrm{T}_x}=\mu X \ts X+ \alpha \w Y\ts \w Y$.
We then define the \emph{apparent pressure} $P_v$ as the pressure $P$ restricted to the horizontal space $H_x$ ; and the \emph{hidden pressure} as $P_h:=P-P_v$.

The fluid will be called perfect if $P_h(Y)=0$. In matrix form, with suitable basis for $T_xH$, $T_xS^1_x$ and $T_x W_x$, and with some abuse of notation for $P_h$, 
$\w G$ can then be written: 
\[
G=\left(\begin{array}{cccc} \mu & \begin{array}{ccc}  0&0&0 \end{array} & e 
\\  \begin{array}{c} 0\\0\\0  \end{array}& P_v & \begin{array}{c} 0\\0\\0 \end{array} & P_h
\\ e &  \begin{array}{ccc} 0&0&0 \end{array}  & \gamma & 0
\\  & P_h & 0 & P_h
\end{array}\right)
\]

This will now be our general model for a fluid.
Then, an analog proof to the above theorem gives :
\begin{greybox}
\begin{theorem} \emph{Equations for the spacetime dynamics of fluids.}

If $D$ is a fluid domain as above, Bianchi identity gives : 
\begin{itemize}
\item Energy Conservation Laws : $\Div (\mu X)=\Div (\mu X_0)=\la X_0,\Div P \ra$ 
\item Electric Charge Conservation Law : $\Div (e X)=\Div (eX_0)=\Div (^eP(Y))= \la Y,\Div P \ra$

\item Motion Equations : 
\begin{itemize}
\item For the fluid : $\mu \wnab_X X=-\Div P-\la X_0,\Div P \ra X$
\item For the apparent fluid : $\mu \wnab_{X_0}X_0 =e.\, ^eF(X_0)-pr_{\mathrm{T}^{\bot}}(\Div P)$
\end{itemize}
\item Maxwell Equations : $\dd F=0$, and $\Div F=e.X_0+1/2|F|_g .Y-\, ^eP(Y)$
\end{itemize}
\end{theorem}
\end{greybox}

\begin{greybox}
Our first theorem above is therefore a special case where the pressure $P=0$.
In this last theorem, the pressure $P$ can be split everywhere into $P=P_v + P_h$, to show physical effects due to the three "classical" dimensions, $P_v$, and those due to the extra "hidden or small" dimensions, $P_h$. This theorem therefore shows that $P_h$, the hidden pressure, is a possible way  to model or encode deviations from standard 4-dimensional General Relativity, or "dark" effects such as dark matter or energy. Indeed, when $P_h=0$, we recover classical equations for a perfect fluid with pressure.
\end{greybox}

Full details of this use of the multi-fibers structure can be found in [3] : S. Collion, M. Vaugon, A New Approach to Kaluza-Klein theory. arXiv.


\end{document}